\documentclass[12pt]{article}
\usepackage{epsfig,amsmath,amsthm}

\newcommand{\x}{{\bf x}}
\newcommand{\z}{{\bf z}}
\newcommand{\y}{{\bf y}}

\newcommand{\n}{{\bf n}}

\newcommand{\e}{{\bf e}}

\newtheorem{Theorem}{Theorem}[section]
\newtheorem{Lemma}{Lemma}[section]

\usepackage {graphics}
\begin{document}

\begin{center}
{\LARGE \bf Homogenization and field concentrations in heterogeneous media}
\end{center}

\begin{center}
{\Large \bf Robert Lipton\\
Department of Mathematics,\\
Louisiana State University, \\
Baton Rouge, LA 70803}
\end{center}

\begin{abstract}
A multi-scale  characterization of the field concentrations inside composite and polycrystalline media is developed.  We focus on gradient fields associated with the intensive quantities given by the
temperature and the electric potential. In the linear regime these quantities are modeled  by the solution of a second order elliptic partial differential equation with oscillatory coefficients.  The characteristic length scale of the heterogeneity relative to the sample size is denoted by $\varepsilon$ and the intensive quantity is denoted by $u^{\varepsilon}$.
Field concentrations are measured using the $L^p$ norm of the gradient field $\Vert\nabla u^\varepsilon\Vert_{\scriptscriptstyle{L^p(D)}}$ for $2\leq p <\infty$. The analysis focuses on the case when $0<\varepsilon \ll 1$.
Explicit lower bounds on $\liminf_{\varepsilon\rightarrow 0}\Vert \nabla u^\varepsilon\Vert_{\scriptscriptstyle{L^p(D)}}$
are developed. 
These bounds provide a way to rigorously assess field concentrations generated by
the microgeometry without having to compute the actual field $u^\varepsilon$.
\end{abstract}




{\bf Abbreviated  title.} Homogenization and field concentrations
\bigskip

\setcounter{equation}{0}
\setcounter{Corollary}{0}
\setcounter{Theorem}{0}
\setcounter{Definition}{0}

\section{\bf Introduction}

The initiation of failure inside heterogeneous media is a multi-scale phenomena. Loads applied at the structural scale are often amplified by the microstructure creating local zones of high field concentration. The local amplification of the applied field creates conditions that are favorable for failure initiation \cite{KellyMac}. This paper focuses on gradient fields associated with the intensive quantities given by the
temperature and the electric potential inside heterogeneous media.  
The local integrability of the
gradient directly correlates with singularity strength which influences the onset of failure such as dielectric breakdown.

In this work it is shown how to assess the $L^p$ 
integrability of the gradient fields in microstructured media by investigating the
multi-scale integrability of suitably defined quantities. 
The analysis is carried out with minimal regularity assumptions on the coefficients describing the local properties inside the heterogeneous media. 
The results are described in terms of the $p^{th}$ order moments of the 
solution of two-scale corrector problems. The quantities are
sensitive to microscopic field concentrations and can become divergent for $p>2$. 
This is in contrast to the well known effective constitutive properties which are based upon local averages and are bounded above independently of the microgeometry.

The results given here are presented in the context of two-scale homogenization 
\cite{allaire}, \cite{Nguetseng}. We consider a bounded domain $\Omega$ in ${\bf R}^n$, $n\geq 2$.
A common microstructure that admits a two-scale description is a simple generalization of a uniformly periodic microstructure and is described as follows.  Consider a partition of
the domain $\Omega$ made up of measurable subsets $\Omega_\ell$, $\ell=1,2,\ldots,K$ such that $\Omega=\cup_{\ell=1}^K\Omega_\ell$. Inside each subdomain $\Omega_\ell$ we place a different periodic microstructure made from $N$ anisotropic heat conductors. This type of microstructure will be referred to
as a piece wise periodic microstructure. Well known  engineering composites that are modeled by  piecewise periodic microstructures include  multi-ply fiber reinforced laminates \cite{gossechristiansen}, \cite{yuan} and \cite{pagano}. 

The thermal conductivity tensor for the piecewise periodic microstructure is described as follows. The indicator function for each of the  subdomains $\Omega_\ell$ is denoted by $\chi_{\scriptscriptstyle{\Omega_\ell}}(\x)$, taking the value $1$ for points in $\Omega_\ell$ and zero outside.  In order to describe the periodic microstructure inside the $\ell^{th}$ subdomain we introduce the unit period cell $Q$. The configuration of the $N$ phases inside $Q$ is described by the indicator functions $\chi_\ell^i(\y)$, $i=1,\ldots,N$ associated with each phase. Here $\chi_\ell^i(\y)=1$ for points inside the $i^{th}$ phase and zero outside. The length scale of the microstructure relative
to the size of the domain $\Omega$ is given by $\varepsilon_k=1/k$, $k=1,2\ldots$. The microstructure is obtained by rescaling the configuration inside the unit period cell.
The indicator function of the $i^{th}$ conductor in the microstructured composite is given by 
\begin{eqnarray}
\chi_i^{\varepsilon_k}(\x)=\chi_i(\x,\x/{\varepsilon_k})=\sum_\ell^K\chi_{\scriptscriptstyle{\Omega_\ell}}(\x)\chi_\ell^i(\x/{\varepsilon_k}).
\label{microindicate}
\end{eqnarray}
The local conductivity tensor $A^{\varepsilon_k}$ has a two-scale structure and is given by
\begin{eqnarray}
A^{\varepsilon_k}(\x)=A(\x,\x/{\varepsilon_k})=\sum_i^N A^i\chi_i(\x,\x/{\varepsilon_k}).
\label{microconduct}
\end{eqnarray}

Other heterogeneous media that are amenable to similar or more general two-scale descriptions include polycrystalline materials such as metals and ceramics. 
We state the general hypotheses under which  the two-scale homogenization theory
applies, see \cite{allaire} and \cite{allairebriane}. 
It is assumed that $A(\x,\y)$ is a matrix defined on $\Omega\times Q$ and there exist positive numbers $\alpha<\beta$ such that
for every vector $\eta$ in ${\bf R}^3$ that
\begin{eqnarray}
\alpha|\eta|^2\leq A(\x,\y)\leq\beta|\eta|^2.
\label{coercivity}
\end{eqnarray}

The conductivity  $A_{ij}(\x,\y)$ is $Q$-periodic in the second variable, such that $A_{ij}(\x,\x/{\varepsilon_k})$ is measurable, satisfies
\begin{eqnarray}
\lim_{{\varepsilon_k}\rightarrow 0}\int_\Omega \left|A_{ij}(\x,\frac{\x}{{\varepsilon_k}})\right|^2\,d\x&=&\int_{\Omega\times Q}\, \left|A_{ij}(\x,\y)\right|^2\,d\x d\y
\label{admiss}
\end{eqnarray}
and for any suitable two scale trial field $\psi(\x,\y)$ that
\begin{eqnarray}
\lim_{{\varepsilon_k}\rightarrow 0}\int_\Omega A_{ij}(\x,\frac{\x}{{\varepsilon_k}})\psi(x,\x/{\varepsilon_k})\,d\x&=&\int_{\Omega\times Q}\, A_{ij}(\x,\y)\psi(\x,\y)\,d\x d\y.
\label{admisstwoscale}
\end{eqnarray}
The convergence given by (\ref{admisstwoscale}) is a  weak convergence and is known as two-scale convergence \cite{allaire} \cite{Nguetseng}. The space of suitable two-scale trials is denoted by $L^2[D;C_{per}(Q)]$.
Here $C_{per}(Q)$ denotes $Q$-periodic continuous functions defined on ${\bf R}^3$
and the space $L^2[D;C_{per}(Q)]$ is the space of functions $h:\Omega\rightarrow C_{per}(Q)$
which are measurable and satisfy $\int_\Omega\Vert h(\x)\Vert^2_{\scriptscriptstyle{C_{per}(Q)}}d\x<\infty$.
The norm $\Vert h(\x)\Vert_{\scriptscriptstyle{C_{per}(Q)}}$ is defined by  $\sup_{\y\in Q}|h(\x,\y)|$. 
In what follows no other regularity hypothesis on the conductivity matrix $A(\x,\y)$ is made.

The temperature field $u^{\varepsilon_k}$ associated with the conductivity tensor field $A^{\varepsilon_k}(\x)=A(\x,\x/{\varepsilon_k})$ is the solution of the
equilibrium equation 
\begin{eqnarray}
-{\rm div}\,\left(A^{\varepsilon_k}(\x) \nabla u^{\varepsilon_k}\right)=f\,\,\,\,{\rm in}\,\,\,\Omega
\label{equilibtwoscale}
\end{eqnarray}
with the boundary conditions given by $u^{\varepsilon_k}\,=\,0$ on $\partial\Omega_D$ and $\n\cdot A^{\varepsilon_k}\nabla u^{\varepsilon_k}=g$ on $\partial\Omega_N$ with $\partial\Omega=\partial\Omega_D\cup\partial\Omega_N$.

In what follows we consider the limit as $\varepsilon_k$ tends to zero. We fix a subdomain $D$ of $\Omega$ and derive lower bounds on 
\begin{eqnarray}
\liminf_{\varepsilon_k\rightarrow 0}\Vert\nabla u^{\varepsilon_k}\Vert_{\scriptscriptstyle{L^p(D)}}.
\label{supinfphase}
\end{eqnarray}

The lower bound is expressed in terms of a two-scale integral that encodes the field amplification properties of the microstructure. It is formulated in terms of the solution of the homogenized problem together with
a local corrector matrix that captures the interaction between the periodic microstructure and the gradients of the homogenized temperature field. The bounds introduced here provide a rigorous way to assess  field concentrations generated by the microgeometry without having to compute the full solution $u^{\varepsilon_k}$.

The lower bound is given in terms of the solutions $w^i(\x,\y)$ to
the local periodic problem. 
For each $\x$ in $\Omega$ the function $w^i(\x,\y)$ is a $Q$ periodic function of the second variable $\y$ and is a solution of
\begin{eqnarray}
{\rm div\,}_{\y}\left(A(\x,\y)(\nabla_\y w^i(\x,\y)+\e^i)\right)=0,
\label{correctorsfuncxt}
\end{eqnarray}
with $\int_Q\,w^i(\x,\y)\,d\y=0$.
The corrector matrix $P(\x,\y)$ is defined by 
\begin{eqnarray}
P_{ij}(\x,\y)=\partial_{y_j}w^i(\x,\y)+e_j^i.
\label{matrixcorrectx}
\end{eqnarray}
The associated effective conductivity tensor $A^E(\x)$ is given by
\begin{eqnarray}
A^E(\x)=\int_Q\, A(\x,\y)P(\x,\y)\,d\y.
\label{effectx}
\end{eqnarray}

The two-scale homogenization theory gives
the following theorem \cite{allaire}.

\begin{Theorem}
{\rm Two-scale Homogenization Theorem}\\
The sequence of solutions $\{u^{\varepsilon_k}\}_{{\varepsilon_k}>0}$ of (\ref{equilibtwoscale})
converges weakly to $u^H(\x)$ in $H^1(\Omega)$ where $u^H$ 
is the solution of the homogenized problem
\begin{eqnarray}
-{\rm div}\,\left( A^E(\x)\nabla u^H(\x)\right)&=&f(\x),\,\,\hbox{in}\,\,\Omega,\nonumber\\
u^H(\x)&=&0 ,\,\,\hbox{on}\,\,\partial\Omega_D,\,\,\hbox{and}\,\,\nonumber\\
 \n\cdot A^E\nabla u^H&=&g,\,\, \hbox{on}\,\, \partial\Omega_N.
\label{homogtwo}
\end{eqnarray}
\label{homogone}
\end{Theorem}

The field concentration functions of
order $p$ are defined by
\begin{eqnarray}
f_p(\x,\nabla u^H(\x))\equiv \left(\int_Q\,|P(\x,\y)\nabla u^H(\x)|^p\,d\y\right)^{1/p}, \,\,\,\,\,2\leq p\leq\infty
\label{macropordx}
\end{eqnarray}
and $f_p(\x,\nabla u^H(\x))\leq f_q(\x,\nabla u^H(\x))$ for $p\leq q$.
It is clear that $f_p$ corresponds to a $p^{th}$ order moment of the corrector matrix (\ref{matrixcorrectx})  and
\begin{eqnarray}
f_\infty(\x,\nabla u^H(\x))&\equiv &\lim_{p\rightarrow\infty}\left(\int_Q\,|P(\x,\y)\nabla u^H(\x)|^p\,d\y\right)^{1/p}.
\label{limsupinftyx}
\end{eqnarray}

\begin{Theorem}
{\rm Lower Bounds on Field Concentrations}\\
For $2 \leq p <\infty$
\begin{eqnarray}
\left(\int_D\,\left(f_p(\x,\nabla u^H(\x))\right)^p\,d\x\right)^{1/p}&\leq& \liminf_{{\varepsilon_k}\rightarrow 0}\Vert\nabla u^{\varepsilon_k}\Vert_{\scriptscriptstyle{L^p(D)}}.
\label{perlowbdx}
\end{eqnarray}
\label{thmi}
\end{Theorem}

For multi-phase conductivity problems with coefficients described by (\ref{microconduct})
the field concentration functions of
order $p$ are defined for each phase and are given by
\begin{eqnarray}
f_p^i(\x,\nabla u^H(\x))\equiv \left(\int_Q\,\chi_i(\x,\y)|P(\x,\y)\nabla u^H(\x)|^p\,d\y\right)^{1/p}, \,\,\,i=1,\ldots,N,\,\,2\leq p\leq\infty
\label{macropordix}
\end{eqnarray}
and $f_p^i(\x,\nabla u^H(\x))\leq f_q^i(\x,\nabla u^H(\x))$ for $p\leq q$.
As before one defines
\begin{eqnarray}
f_\infty^i(\x,\nabla u^H(\x))&\equiv &\lim_{p\rightarrow\infty}\left(\int_Q\,\chi_i(\x,\y)|P(\x,\y)\nabla u^H(\x)|^p\,d\y\right)^{1/p}.
\label{limsupinftyix}
\end{eqnarray}
For this case lower bounds on 
\begin{eqnarray}
\liminf_{\varepsilon_k\rightarrow 0}\Vert\chi_i^{\varepsilon_k}\nabla u^{\varepsilon_k}\Vert_{\scriptscriptstyle{L^p(D)}}
\label{supinfphasei}
\end{eqnarray}
are given by the following theorem.

\begin{Theorem}
{\rm Lower Bounds for Multi-phase Composites}\\
For $2 \leq p <\infty$
\begin{eqnarray}
\left(\int_D\,\left(f^i_p(\x,\nabla u^H(\x))\right)^p\,d\x\right)^{1/p}&\leq& \liminf_{{\varepsilon_k}\rightarrow 0}\Vert\chi_i^{\varepsilon_k}\nabla u^{\varepsilon_k}\Vert_{\scriptscriptstyle{L^p(D)}}.
\label{perlowbdxi}
\end{eqnarray}
\label{thmmi}
\end{Theorem}

The bounds can be applied to develop a Chebyshev Inequality for the distribution functions associated with
the sequence $\{\chi_i^{\varepsilon_k}|\nabla u^{\varepsilon_k}|\}_{{\varepsilon_k}>0}$.
Here the distribution function $\lambda_i^{\varepsilon_k}(D,t)$ gives the measure of the set inside $D$ where $\chi_i^{\varepsilon_k}|\nabla u^{\varepsilon_k}|>t$.

Arguing as in  Proposition 2.1 of \cite{lipsiap} and combining with (\ref{perlowbdxi}) gives the following
\begin{Theorem}
{\rm Homogenized Chebyshev Inequality}\\
\begin{eqnarray}
\limsup_{\varepsilon_k\rightarrow 0}\lambda_i^{\varepsilon_k}(D,t)&\leq t^{-p}\left(\int_D\,\left(f_p^i(\x,\nabla u^H(\x))\right)^p\,d\x\right)&\leq t^{-p} \liminf_{{\varepsilon_k}\rightarrow 0}\Vert\chi_i^{\varepsilon_k}\nabla u^{\varepsilon_k}\Vert^p_{\scriptscriptstyle{L^p(D)}}.\nonumber\\
\label{perChebyshevxi}
\end{eqnarray}
\label{Chebyshevtwoscale}
\end{Theorem}

It is pointed out that Theorems \ref{thmi}  and  \ref{thmmi} are obtained using the minimum regularity assumptions on the coefficients $A^{\varepsilon_k}$. Because of this the hypotheses of (Theorem 2.6, \cite{allaire}) do not apply and one can not take advantage of the strong convergence given in that Theorem. Instead the theorems are proved using a perturbation approach
introduced in \cite{liproy}, \cite{lipquart}, see Section 2.

The lower bounds are sensitive to the presence of singularities generated by the microstructure.
To illustrate this we consider a microstructure made from a periodic distribution of uniaxial crystallites embedded in an isotropic matrix of unit conductivity. The period cell for the composite is illustrated in Figure 1. Each crystallite occupies a sphere and has conductivity $\lambda_1$ in the radial direction and $\lambda_2$ in the tangential direction. The dispersion of the $N$ crystallites is specified by $\cup_\ell^N B(\y^\ell,r_\ell)$ where $B(\y^\ell,r_\ell)$ denotes the $\ell^{th}$ sphere centered
at $\y^\ell$ with radius $r_\ell$. Each crystallite has a conductivity tensor given by
\begin{eqnarray}
A(\y)=\lambda_1\n\otimes\n+\lambda_2(I-\n\otimes\n),
\label{loccondct}
\end{eqnarray}
where
$\n=(\y-\y^\ell)/|\y-\y^\ell|$ for $\y$ in $B(\y^\ell,r_\ell)$ and $I$ is the $3\times 3$ identity. 
Outside the crystallites we set $A(\y)=I$. It is supposed that the aggregate of crystallites occupy an area fraction $0<\theta<1$ of the unit period cell. It is noted that the conductivity inside each crystallite is precisely the one employed in the Schulgasser sphere assemblage \cite{schul}. 

When a constant gradient field
is applied to a single  isolated crystallite and $\lambda_1>\lambda_2$ 
the crystallite exhibits a gradient field singularity at its center. 
In what follows we use the lower bound (\ref{perlowbdx}) to show
how this local information effects the integrability of the sequence $\{\nabla u^{\varepsilon_k}\}_{\varepsilon_k>0}$.
We form $A^{\varepsilon_k}=A(\x/\varepsilon_k)$ and consider solutions $u^{\varepsilon_k}$ of 
(\ref{equilibtwoscale}). To fix ideas we choose $f$ to be in $L^r(\Omega)$ for $r>3$ and $g$ to be in $L^2(\partial\Omega_N)$.
In what follows $\lambda_2$ is restricted to lie in the interval $1/2<\lambda_2<1$ and  $\lambda_1=1/(2\lambda_2-1)$. 
For this choice it is shown in Section 3 that the homogenized temperature field $u^H$ is the solution of (\ref{homogtwo}) with $A^E=I$.

For $D$ compactly contained in $\Omega$ it follows from
the $L^p$ theory \cite{meyers}, that $\Vert \nabla u^H\Vert_{L^p(D)}<\infty$ for every $1\leq p<\infty$.
On the other hand calculation and application of Theorem \ref{thmi} shows that
\begin{eqnarray}
LB(p)\times\Vert\nabla u^H\Vert_{\scriptscriptstyle{L^p(D)}}&\leq& \liminf_{{\varepsilon_k}\rightarrow 0}\Vert\nabla u^{\varepsilon_k}\Vert_{\scriptscriptstyle{L^p(D)}},
\label{lowest}
\end{eqnarray}
where
\begin{eqnarray}
LB(p)=\left\{\begin{array}{l c}
\frac{3p\theta(2\lambda_2-1)}{2(1-\lambda_2)(\frac{3}{2(1-\lambda_2)}-p)}+(1-\theta),&\hbox{for},\,\,p<\frac{3}{2(1-\lambda_2)},\\
+\infty,&\hbox{for},\,\,p\geq\frac{3}{2(1-\lambda_2)}
\end{array}.\right.
\label{lowdivergence}
\end{eqnarray}
For a fixed choice of $\lambda_2$ the value $p_c=\frac{3}{2(1-\lambda_2)}$ satisfies $3< p_c<+\infty$ and 
\begin{eqnarray}
\liminf_{{\varepsilon_k}\rightarrow 0}\Vert\nabla u^{\varepsilon_k}\Vert_{\scriptscriptstyle{L^p(D)}}&=&+\infty,\,\,\hbox{for}\,\,p\geq p_c.
\label{lowestdiverge}
\end{eqnarray}
This is in stark contrast to the $L^p$ integrability of the the gradient of the homogenized solution
which holds for any $p<+\infty$.
It is clear for this example that the information carried by the homogenized problem is not adequate and misses the singular behavior exhibited by the sequence $\{\nabla u^{\varepsilon_k}\}_{\varepsilon_k>0}$. This example shows that failure initiation criteria based solely upon the solution of the homogenized equations will be optimistic.
The inequalities given above are established in Section 3.

The maximum integrability exponent for the gradient of the solution of the local problem (\ref{correctorsfuncxt}) is referred to as
the threshold exponent for the composite. The threshold exponent is introduced in the work of  Milton
\cite{milton2} and measures the worst singularity of the gradient field. The threshold 
exponent is precisely $p_c$ for the local problem considered here and corresponds to the divergence in the lower bound
for $p\geq p_c$.

\begin{figure}[htb]
\centerline{\epsfig{file=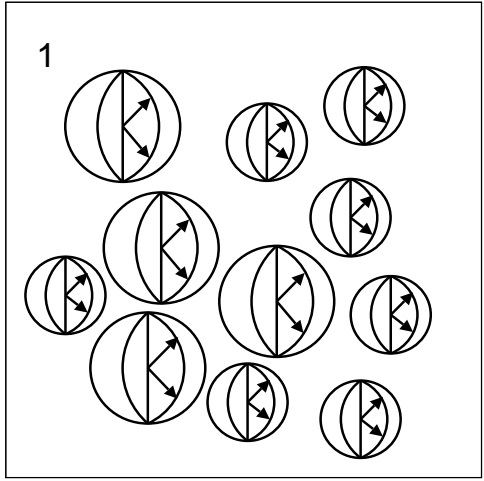,height=2.5in,width=2.5in}}
\caption{Unit period cell with Schulgasser crystallites embedded  inside a material with unit thermal conductivity.} \label{108}
\end{figure}

Next we consider an example for which 
the sequence  $\{\nabla u^{\varepsilon_k}\}_{\varepsilon_k>0}$ is uniformly bounded in $L^p$ 
for some class of coefficients and right hand sides $f$. For this case we show that
the lower bound given in Theorem \ref{thmi}
is attained. In this example we make use of the a priori estimates for $\{\nabla u^{\varepsilon_k}\}_{\varepsilon_k>0}$
developed in the Theorem 4 of Avellaneda and Lin \cite{AL}. Let $\Omega$ be a $C^{1,\alpha}$ domain $(0<\alpha\leq 1)$ and suppose for $0<\gamma\leq 1$, $C>0$, that $A(\y)\in C^\gamma({\bf R}^n)$
and $\Vert A(\y)\Vert_{C^\gamma({\bf R}^n)}\leq C$. Then we choose $A^{\varepsilon_k}=A(\x/\varepsilon_k)$.
For $\delta>0$ suppose $2\leq q\leq n+\delta$ and $f\in L^q$ and set  $1/\hat{q}=1/q-1/(n+\delta)$.
Given these choices we consider the $W^{1,2}_0(\Omega)$ solutions $u^{\varepsilon_k}$ 
of
\begin{eqnarray}
-{\rm div}\,\left(A^{\varepsilon_k}(\x) \nabla u^{\varepsilon_k}\right)=f\,\,\,\,{\rm in}\,\,\,\Omega.
\label{equilibtwoscalex}
\end{eqnarray}
It is shown in Section 4 that (\ref{perlowbdx}) holds with equality
for every $p$ such that $p<\hat{q}$.
In fact it is seen more generally that for $p<\hat{q}$ and any Caratheodory function $\psi:D\times{\bf R}^n\rightarrow{\bf R}$ satisfying
\begin{eqnarray}
|\psi(\x,\eta)|\leq |\eta|^p,\,\,\hbox{for a.e.  }\x\in D \,\,\hbox{and  }\eta\in {\bf R}^3 ,
\label{equality}
\end{eqnarray}
that
\begin{eqnarray}
\lim_{\varepsilon_k\rightarrow 0}\int_D\,\psi(\x,\nabla u^{\varepsilon_k}(\x))\,d\x&=&\int_D\int_Q\,\psi(\x,P(\y)\nabla u^H(\x))\,d\y d\x.
\label{caratheodory}
\end{eqnarray}
This is established in Section 4.

It is anticipated that there are several classes of conductivity coefficients and right hand sides $f$ for
which the lower bounds are attained. 
In this direction we point out the the recent higher regularity results given in 
\cite{bv}, \cite{cp}, \cite{peralguterez}, \cite{lv}, \cite{ln} and \cite{schweizer}.

We conclude noting that the analogues of the field concentration functions (\ref{macropordx}) and (\ref{macropordix}) have appeared earlier in the contexts of G-convergence and random media, see \cite{liproy} and \cite{lipsiap}. In those treatments
they are shown to provide upper bounds for the distribution function of the local stress and electric field for G-convergent sequences of elasticity tensors and random dielectric tensors.

\setcounter{equation}{0}
\setcounter{Corollary}{0}
\setcounter{Theorem}{0}
\setcounter{Definition}{0}

\section{\bf Derivation of the lower bounds}

We recall the weak formulation of the ${\varepsilon_k}>0$
problem given by (\ref{equilibtwoscale}).
Let $V$ denote the closure in  $H^1(\Omega)$ of all smooth functions that vanish on $\partial\Omega_D$. We suppose that $f$ is in $L^2(\Omega)$ and $g$ belongs to $L^2(\partial\Omega_N)$.  The function $u^{\varepsilon_k}$ belonging to $V$  is the solution of the weak formulation of the boundary value problem given by
\begin{eqnarray}
\int_\Omega\,A(\x,\x/{\varepsilon_k})\nabla u^{\varepsilon_k}\cdot\nabla \varphi\,d\x&=&\int_\Omega\,f\varphi\,d\x+\int_{\partial\Omega_N}\,g\varphi\,ds,
\label{weakform}
\end{eqnarray}
for every $\varphi$ in $V$. Here $ds$ is an element of surface area.

In order to express the two-scale weak formulation of (\ref{homogtwo}) we introduce the following 
function spaces. The space of square integrable $Q$-periodic mean zero functions with square integrable derivatives
is denoted by $H^1_{per}(Q)/{\bf R}$. The norm of an element $v$ in this space is denoted by $\Vert v \Vert_{\scriptscriptstyle{H^1_{per}(Q)/{\bf R}}}$.
The space of measurable functions $h$ from $\Omega$ to $H^1_{per}(Q)/{\bf R}$ for which $\int_\Omega\,\Vert h(\x)\Vert^2_{\scriptscriptstyle{H^1_{per}(Q)/{\bf R}}}\,d\x<\infty$ is denoted by $L^2[\Omega;H^1_{per}(Q)/{\bf R}]$. This function space was introduced for the description of the two-scale homogenized problem in \cite{Nguetseng}. The weak formulation of the two-scale homogenized problem (\ref{homogtwo})
is given by the  unfolded variational principle \cite{allaire}, \cite{cioran}, \cite{lukkasenguetsengwall}.

\begin{Theorem}
{\rm Unfolded Variational Principle}\\
The pair $(u^H,u_1)$ is the unique solution in $V\times L^2[\Omega;H^1_{per}(Q)/{\bf R}]$ of 
\begin{eqnarray}
&&\int_\Omega\int_Q\,A(\x,\y)(\nabla u^H(\x) + \nabla_\y u_1(\x,\y))\cdot(\nabla \varphi(\x)+\nabla_\y \varphi_1(\x,\y))\,d\y\,d\x\nonumber\\
&&=\int_\Omega\,f\varphi\,d\x+\int_{\partial\Omega_N}\,g\varphi\,ds,
\label{homogweakform}
\end{eqnarray}
for every $(\varphi,\varphi_1)$ in $V\times L^2[\Omega;H^1_{per}(Q)/{\bf R}]$.  Moreover
\begin{eqnarray}
\nabla u^H+\nabla_{\y}u_1(\x,\y)=P(\x,\y)\nabla u^H(\x).
\label{correctx}
\end{eqnarray}
\label{twohomog}
\end{Theorem}

In order to establish Theorems \ref{thmi} and \ref{thmmi} we recall the function spaces used in the description of two-scale convergence \cite{lukkasenguetsengwall}.
The space $C_{per}(Q)$ denotes $Q$-periodic continuous functions defined on ${\bf R}^3$.
For $1\leq r <\infty$ the space $L^r[D;C_{per}(Q)]$ is the space of functions $h:D \rightarrow C_{per}(Q)$
which are measurable and satisfy $\int_D \Vert h(\x)\Vert^r_{\scriptscriptstyle{C_{per}(Q)}}d\x<\infty$.
Here $\Vert h(\x)\Vert_{\scriptscriptstyle{C_{per}(Q)}}=\sup_{\y\in Q}|h(\x,\y)|$. The intersection of the spaces $L^\infty(D \times Q)$ and $L^r[D;C_{per}(Q)]$ is denoted by by $V^r$. For $1<r<\infty$ we introduce $1<r'<\infty$ such that $\frac{1}{r}+\frac{1}{r'}=1$.
We establish Theorems \ref{thmi} and \ref{thmmi} with the aid of the following Lemmas.

\begin{Lemma}
{\rm Localization Lemma}\\
Fix a domain of interest $D$ inside $\Omega$. Let $q(\x,\y)$ be any test function in $V^{r}$ then one can pass to the limit ${\varepsilon_k}\rightarrow 0$ in the sequence
of solutions $\{u^{\varepsilon_k}\}_{\varepsilon_k >0}$ of (\ref{equilibtwoscale}) to obtain:
\begin{eqnarray}
\lim_{{\varepsilon_k}\rightarrow 0}\int_D\,q(\x,\x/{\varepsilon_k})\,|
\nabla u^{\varepsilon_k}|^2\,d\x
&=& \int_D\,\int_Q\,q(\x,\y)\,|P(\x,\y)\nabla u^H(\x)|^2\,d\y\,d\x.
\label{limloc}
\end{eqnarray}
\label{locallemma}
\end{Lemma}
For multi-phase composites with coefficients described by (\ref{microconduct}) we restrict attention inside each phase and state the following lemma.
\begin{Lemma}
{\rm Localization Lemma in Multi-phase Composites}\\
Let $q(\x,\y)$
be any test function in $V^{r}$ then one can pass to the limit ${\varepsilon_k}\rightarrow 0$ in the sequence
of solutions $\{u^{\varepsilon_k}\}_{{\varepsilon_k}>0}$ of (\ref{equilibtwoscale}) to obtain:
\begin{eqnarray}
&&\lim_{{\varepsilon_k}\rightarrow 0}\int_D\,q(\x,\x/{\varepsilon_k})\,\chi_i^{\varepsilon_k}(\x)\,|
\nabla u^{\varepsilon_k}|^2\,d\x \nonumber\\
&&= \int_D\,\int_Q\,q(\x,\y)\,\chi_i(\x,\y)\,|P(\x,\y)\nabla u^H(\x)|^2\,d\y\,d\x.
\label{limloci}
\end{eqnarray}
\label{locallemmai}
\end{Lemma}
\noindent The proofs of Lemmas \ref{locallemma} and \ref{locallemmai} are given at the end of this section.

To illustrate the ideas we use Lemma \ref{locallemmai}  to establish  Theorem \ref{thmmi} noting that Theorem \ref{thmi} follows from Lemma \ref{locallemma}  in the same way. 

{\bf Proof of Theorem \ref{thmmi}.}
For each ${\varepsilon_k}>0$ we apply H\"older's inequality to the left side of (\ref{limloci}) to obtain
\begin{eqnarray}
&&\int_D\,\int_Q\,q(\x,\y)\,\chi_i(\x,\y)\,|P(\x,\y)\nabla u^H(\x)|^2\,d\y\,d\x\nonumber\\
&&\leq\lim_{{\varepsilon_k}\rightarrow 0}\left(\int_D\,\vert q(\x,\x/{\varepsilon_k})\vert^{r}\,dx\right)^{1/r}\liminf_{{\varepsilon_k}\rightarrow 0}\left(\int_D\,\chi_i^{\varepsilon_k}(\x)\,|
\nabla u^{\varepsilon_k}|^{2r'}\,d\x\right)^{1/r'}. 
\label{limlociup}
\end{eqnarray}
Noting \cite{lukkasenguetsengwall} that
\begin{eqnarray}
\lim_{{\varepsilon_k}\rightarrow 0}\left(\int_D\,|q(\x,\x/{\varepsilon_k})|^r\,dx\right)^{1/r}
=& \left(\int_D\,\int_Q\,|q(\x,\y)|^r\,d\y\,d\x\right)^{1/r}&\equiv \Vert q(\x,\y)\Vert_{\scriptscriptstyle{L^{r}(D\times Q)}}
\label{limdecouple}
\end{eqnarray}
we obtain
\begin{eqnarray}
\frac{\int_D\,\int_Q\,q(\x,\y)\,\chi_i(\x,\y)\,|P(\x,\y)\nabla u^H(\x)|^2\,d\y\,d\x}{\Vert q(\x,\y)\Vert_{\scriptscriptstyle{L^{r}(D\times Q)}}}
&\leq&\liminf_{{\varepsilon_k}\rightarrow 0}\left(\int_D\chi_i^{\varepsilon_k}(\x)|
\nabla u^{\varepsilon_k}|^{2r'}d\x\right)^{1/r'}.
\label{preestimate}
\end{eqnarray}
Since $V^r$ is dense in $L^r(D\times Q)$ we take the 
supremum of the left hand side of (\ref{preestimate}) over 
$V^r$ to find that
\begin{eqnarray}
\left(\int_D\,\int_Q\,\chi_i(\x,\y)\,|P(\x,\y)\nabla u^H(\x)|^{2r'}\,d\y\,d\x\right)^{1/r'}
&\leq&\liminf_{{\varepsilon_k}\rightarrow 0}\left(\int_D\chi_i^{\varepsilon_k}(\x)|
\nabla u^{\varepsilon_k}|^{2r'}d\x\right)^{1/r'}.
\label{eestimate}
\end{eqnarray}
Theorem \ref{thmmi} follows for $2<p<\infty$ upon taking the square root on both sides of (\ref{eestimate}).
The case $p=2$ follows immediately upon choosing $q(\x,\y)=1$ in  Lemma \ref{locallemmai}.

We conclude by providing the proof of Lemma \ref{locallemmai} and note that the proof of Lemma \ref{locallemma}
is identical. 

{\bf Proof of Lemma \ref{locallemmai}.}
The indicator function of the set of interest $D$ is denoted  by $\chi_D(\x)$. We choose a test function $q(\x,\y)$ in $V^r$ and set $p(\x,\y)=\chi_D(\x)\chi_i(\x,\y) q(\x,\y)$. For $\delta\beta >0$ we form the perturbed conductivity tensor $\tilde{A}_{ij}(\x,\y)=A_{ij}(\x,\y)+\delta\beta p(\x,\y)\delta_{ij}$. We choose $\delta\beta$ sufficiently small so that  $\tilde{A}(\x,\y)$ satisfies (\ref{coercivity}). By construction $\tilde{A}(\x,\x/{\varepsilon_k})$ is measurable and satisfies (\ref{admiss}) and (\ref{admisstwoscale}).
Consider the associated solution $\tilde{u}^{\varepsilon_k}$ in $V$  of the weak formulation of the boundary value problem given by
\begin{eqnarray}
\int_\Omega\,\tilde{A}(\x,\x/{\varepsilon_k})\nabla \tilde{u}^{\varepsilon_k}\cdot\nabla \varphi\,d\x&=&\int_\Omega\,f\varphi\,d\x+\int_{\partial\Omega_N}\,g\varphi\,ds,\, \,\,\,\hbox{for every $\varphi$ in $V$.}
\label{weakformpurt}
\end{eqnarray}
Set $\tilde{u}^{\varepsilon_k}=u^{\varepsilon_k}+\delta u^{\varepsilon_k}$ and subtraction of (\ref{weakform}) from
(\ref{weakformpurt}) gives
\begin{eqnarray}
\int_\Omega\,\tilde{A}(\x,\x/{\varepsilon_k})\nabla \delta {u}^{\varepsilon_k}\cdot\nabla \varphi\,d\x+
\int_\Omega\,\delta\beta\,p(\x,\x/{\varepsilon_k})\nabla {u}^{\varepsilon_k}\cdot\nabla \varphi\,d\x&=& 0.
\label{weakformvar}
\end{eqnarray}
Choosing $\varphi=u^{\varepsilon_k}$ in (\ref{weakformvar}) and application of the identity
\begin{eqnarray}
\int_\Omega\,A(\x,\x/{\varepsilon_k})\nabla u^{\varepsilon_k}\cdot\nabla \delta u^{\varepsilon_k}\,d\x&=&\int_\Omega\,f\delta u^{\varepsilon_k}\,d\x+\int_{\partial\Omega_N}\,g\delta u^{\varepsilon_k}\,ds,
\label{crossterm}
\end{eqnarray}
gives
\begin{eqnarray}
\delta\beta\times\int_\Omega\,p(\x,\x/{\varepsilon_k})|\nabla u^{\varepsilon_k}|^2\,d\x + T^{\varepsilon_k} &=&-\int_\Omega\,f\delta u^{\varepsilon_k}\,d\x-\int_{\partial\Omega_N}\,g\delta u^{\varepsilon_k}\,ds,
\label{epsdelta}
\end{eqnarray}
where
\begin{eqnarray}
T^{\varepsilon_k}&=&\delta\beta\times\int_\Omega\,p(\x,\x/{\varepsilon_k})(\nabla \delta u^{\varepsilon_k})\cdot\nabla u^{\varepsilon_k}\,d\x.
\end{eqnarray}
Next set $\varphi=\delta u^{\varepsilon_k}$ in (\ref{weakformvar})
and it follows from Cauchy's inequality and (\ref{coercivity}) that 
\begin{eqnarray}
\Vert\nabla \delta u^{\varepsilon_k}\Vert_{\scriptscriptstyle{L^2(\Omega)}}\leq C\delta\beta,
\label{std}
\end{eqnarray}
where here and throughout $C$ denotes a generic a constant independent of ${\varepsilon_k}$. From this it is evident that 
\begin{eqnarray}
|T^{\varepsilon_k}|<C\delta\beta^2.
\label{betaest}
\end{eqnarray}
Next we pass to the ${\varepsilon_k}\rightarrow 0$ limit and apply Theorems \ref{homogone} and \ref{twohomog}  to find that the sequence $\{\tilde{u}^{\varepsilon_k}\}_{{\varepsilon_k}>0}$ converges weakly in $H^1(\Omega)$ to
$\tilde{u}^H$, where $(\tilde{u}^H,\tilde{u}_1)$ is the solution in $V\times L^2[\Omega;H^1_{per}(Q)/{\bf R}]$ of 
\begin{eqnarray}
&&\int_\Omega\int_Q\,\tilde{A}(\x,\y)(\nabla \tilde{u}^H(\x) + \nabla_\y \tilde{u}_1(\x,\y))\cdot(\nabla \varphi(\x)+\nabla_\y \varphi_1(\x,\y))\,d\y\,d\x\nonumber\\
&&=\int_\Omega\,f\varphi\,d\x+\int_{\partial\Omega_N}\,g\varphi\,ds,
\label{homogweakformpurt}
\end{eqnarray}
for every $(\varphi,\varphi_1)$ in $V\times L^2[\Omega;H^1_{per}(Q)/{\bf R}]$.
Set $\tilde{u}^H-u^H=\delta u^H$, $\tilde{u}_1-u_1=\delta u_1$ and subtraction of (\ref{homogweakform}) from (\ref{homogweakformpurt}) gives
\begin{eqnarray}
&&\int_\Omega\int_Q\,\tilde{A}(\x,\y)(\nabla \delta{u}^H(\x) + \nabla_\y \delta{u}_1(\x,\y))\cdot(\nabla \varphi(\x)+\nabla_\y \varphi_1(\x,\y))\,d\y\,d\x\nonumber\\
&&+\int_\Omega\int_Q\,\delta\beta p(\x,\y)(\nabla {u}^H(\x) + \nabla_\y {u}_1(\x,\y))\cdot(\nabla \varphi(\x)+\nabla_\y \varphi_1(\x,\y))\,d\y\,d\x = 0.
\label{homogweakformvar}
\end{eqnarray}
Choosing $(\varphi,\varphi_1)=({u}^H,{u}_1)$ in (\ref{homogweakformvar}) together
with the identity 
\begin{eqnarray}
&&\int_\Omega\int_Q\,{A}(\x,\y)(\nabla {u}^H(\x) + \nabla_\y {u}_1(\x,\y))\cdot(\nabla\delta u^H(\x)+\nabla_{\y}\delta u_1(\x,\y))\,d\y\,d\x\nonumber\\
&&=\int_\Omega\,f\delta {u}^H\,d\x+\int_{\partial\Omega_N}\,g\delta {u}^H\,ds.
\label{homogcrossterm}
\end{eqnarray}
gives
\begin{eqnarray}
&&\delta\beta\times\int_\Omega\int_Q\,p(\x,\y)|P(\x,\y)\nabla{u}^H(\x)|^2\,d\y\,d\x + \tilde{T}\nonumber\\
&&=-\int_\Omega\,f\delta u^H\,d\x-\int_{\partial\Omega_N}\,g\delta u^H\,ds,
\label{homogweakorders}
\end{eqnarray}
where
\begin{eqnarray}
\tilde{T} &=&\delta\beta\times\int_\Omega\int_Q\,p(\x,\y)(\nabla \delta u^H+\nabla_{\y}\delta u_1(\x,\y))\cdot(\nabla u^H+\nabla_{\y}u_1(\x,\y))\,d\x.
\label{highorderhomog}
\end{eqnarray}
Next set $(\varphi,\varphi_1)=(\delta{u}^H,\delta{u}_1)$ in (\ref{homogweakformvar})
and it follows from Cauchy's inequality and (\ref{coercivity}) that 
\begin{eqnarray}
\Vert\nabla \delta u^H+\nabla_{y}\delta u_1\Vert_{\scriptscriptstyle{L^2(\Omega\times Q)}}\leq C\delta\beta
\label{stdh}
\end{eqnarray}
and it follows easily that 
\begin{eqnarray}
|\tilde{T}|<C\delta\beta^2.
\label{betaesth}
\end{eqnarray}
Taking the ${\varepsilon_k}\rightarrow 0$ limit in (\ref{epsdelta}) noting that $\lim_{{\varepsilon_k}\rightarrow 0}\delta u^{\varepsilon_k}=\delta u^H$ (weakly in $H^1(\Omega)$) and recalling (\ref{betaest}) gives
\begin{eqnarray}
\delta\beta\times\lim_{{\varepsilon_k}\rightarrow 0}\int_\Omega\,p(\x,\x/{\varepsilon_k})|\nabla u^{\varepsilon_k}|^2\,d\x + O(\delta\beta^2) &=&-\int_\Omega\,f\delta u^H\,d\x-\int_{\partial\Omega_N}\,g\delta u^H ds.
\label{epsdeltal}
\end{eqnarray}
Lemma \ref{locallemmai} now follows immediately from (\ref{homogweakorders}),  (\ref{betaesth}) and (\ref{epsdeltal}) and identifying like powers of $\delta\beta$.

\setcounter{equation}{0}
\setcounter{Corollary}{0}
\setcounter{Theorem}{0}
\setcounter{Definition}{0}

\section{\bf Explicit lower bounds for  aggregates of Schulgasser crystallites}

In this section we derive the lower bound (\ref{lowdivergence}) for the microstructure consisting of  Schulgasser
crystallites embedded within a homogeneous matrix with unit thermal conductivity. 
The temperature field inside the unit period cell $\Phi^i(\y)=w^i(\y)+\y_i$ is the solution
of the local problem
\begin{eqnarray}
{\rm div\,}_{\y}\left(A(\y)(\nabla_\y w^i(\y)+\e^i)\right)=0,
\label{correctorsfuncxtc}
\end{eqnarray}
with $w^i$ Q-periodic and $\int_Q\,w^i(\y)\,d\y=0$.
For this microstructure $A(\y)$ is given by (\ref{loccondct}) for $\y$ in $B(\y^\ell,r_\ell)$ and $A(\y)=I$
outside. Here we restrict $\lambda_2$ to the interval $(1/2,1)$ and choose $\lambda_1$ so that $\lambda_1=1/(2\lambda_2-1)$. A calculation shows that the solution $\Phi^i(\y)$ is given by
\begin{eqnarray}
\Phi^i=\left\{\begin{array}{c l}
\y_i,&\y\in Q\setminus\cup_{\ell=1}^N B(\y^\ell,r_\ell),\\
r_\ell^{1-\alpha}|\y-\y^\ell|^{\alpha-1}(\y_i-\y^\ell_i)+\y_i^\ell,&\y\in B(\y^\ell,r_\ell)
\end{array}\right.,
\label{soln}
\end{eqnarray}
where $\alpha=2\lambda_2-1$.
The corrector matrix $P(\y)$ is given by
\begin{eqnarray}
P(\y)=\left\{\begin{array}{c l}
I,&\y\in Q\setminus\cup_{\ell=1}^N B(\y^\ell,r_\ell),\\
r_\ell^{1-\alpha}|\y-\y^\ell|^{\alpha-1}(I+(\alpha-1)\n\otimes\n),&\y\in B(\y^\ell,r_\ell)
\end{array}\right.,
\label{solncorrect}
\end{eqnarray}
where $\n=(\y-\y^\ell)/|\y-\y^\ell|$ for $\y\in B(\y^\ell,r_\ell)$.
A direct calculation shows that 
\begin{eqnarray}
A^E=\int_Q\,A(\y)P(\y)\,d\y=I.
\label{effect}
\end{eqnarray}
Next we provide the lower bound for $\int_\Omega\int_Q\,|P(\y)\nabla u^H(\x)|^p\,d\y\,d\x$.
Note for any $\eta$ in ${\bf R}^3$ that $P^T(\y)P(\y)\eta\cdot\eta=|P(\y)\eta|^2$ and the smallest
eigenvalue $\lambda(\y)$ of $P^T(\y)P(\y)$ delivers the lower bound $\lambda(\y)|\eta|^2\leq|P(\y)\eta|^2$
and 
\begin{equation}
\int_\Omega\,\int_Q\,\lambda(\y)^{p/2}|\nabla u^H(\x)|^p\,d\y\,d\x\leq\int_\Omega\,\int_Q\,|P(\y)\nabla u^H(\x)|^p\,d\y\,d\x. 
\label{lowbddd}
\end{equation}
Calculation shows that 
\begin{eqnarray}
\lambda(\y)=\alpha^2 r_\ell^{2(1-\alpha)} |\y-\y^\ell|^{2(\alpha-1)}
\label{lambda}
\end{eqnarray}
for $\y\in B(\y^\ell,r_\ell)$ and $\lambda(\y)=1$ for $\y\in Q\setminus \cup_\ell^N B(\y^\ell,r_\ell)$.
The lower bound (\ref{lowdivergence}) follows upon substitution of (\ref{lambda}) into (\ref{lowbddd}).

\setcounter{equation}{0}
\setcounter{Corollary}{0}
\setcounter{Theorem}{0}
\setcounter{Definition}{0}

\section{\bf Optimality of the lower bounds}

Conditions are presented on  $f$ and $A(\y)$ for which the lower bound 
(\ref{perlowbdx}) is attained
for a range of exponents $2< p<\hat{q}$. We suppose as in Avellaneda and Lin \cite{AL} that $\Omega$ is a $C^{1,\alpha}$ domain $(0<\alpha\leq 1)$ and suppose for $0<\gamma\leq 1$, $0<C$, that $A(\y)\in C^\gamma({\bf R}^n)$
and $\Vert A(\y)\Vert_{C^\gamma({\bf R}^n)}\leq C$. We set $A^{\varepsilon_k}=A(\x/\varepsilon_k)$.
For $\delta>0$ suppose $2\leq q\leq n+\delta$ and $f\in L^q$ and set  $1/\hat{q}=1/q-1/(n+\delta)$.
Given these choices we consider the $W^{1,2}_0(\Omega)$ solutions $u^{\varepsilon_k}$ 
of
\begin{eqnarray}
-{\rm div}\,\left(A^{\varepsilon_k}(\x) \nabla u^{\varepsilon_k}\right)=f\,\,\,\,{\rm in}\,\,\,\Omega.
\label{equilibtwoscalex}
\end{eqnarray}
Theorem 4 of \cite{AL} shows that there exists a constant independent of $\varepsilon_k$ for which
\begin{eqnarray}
\Vert\nabla u^{\varepsilon_k}\Vert_{L^{\hat{q}}(\Omega)}\leq C\Vert f\Vert_{L^q(\Omega)}
\label{upperboundek}
\end{eqnarray}
holds for every $\varepsilon_k>0$.
Subject to these hypotheses it will be shown that the lower bound (\ref{perlowbdx}) is attained for $p<\hat{q}$.

Passing to a subsequence if necessary we start by considering the Young measure $\nu$ associated with the sequence $\{P(\x/\varepsilon_k)\nabla u^H(\x)\}_{\varepsilon_k >0}$. Here $\nu$ is represented by a family of probability measures $\nu=\{\nu_\x\}_{\x\in\Omega}$ depending measurably on $\x$.  We denote the set of continuous functions $\varphi$ defined on ${\bf R}^n$
such that $\lim_{\eta\rightarrow\infty}\varphi(\eta)=0$ by $C_0({\bf R}^n)$.  Elementary arguments show that
\begin{eqnarray}
<\nu_\x,\varphi>=\int_{{\bf R}^n}\varphi(\eta)d\nu_\x(\eta)=\int_Q\,\varphi(P(\z)\nabla u^H(\x))d\z,\,\,\,\hbox{a.e.}\,\,\x\in\Omega
\label{identify}
\end{eqnarray}
for every $\varphi$ in $C_0({\bf R}^n)$.
From corrector theory \cite{murattartar} there exists an exponent $r\geq 1$ for which one has the strong convergence
\begin{eqnarray}
\lim_{\varepsilon_k\rightarrow 0}\Vert\nabla u^{\varepsilon_k}-P(\x/\varepsilon_k)\nabla u^H\Vert_{L^r(\Omega)}=0.
\label{stongconv}
\end{eqnarray}
The strong convergence (\ref{stongconv}) shows that both sequences $\{\nabla u^{\varepsilon_k}\}_{\varepsilon_k >0}$ and $\{P(\x/\varepsilon_k)\nabla u^H(\x)\}_{\varepsilon_k >0}$ share the same Young measure see for example Lemma 6.3 of \cite{pedregal}.
From (\ref{upperboundek}) it follows on passage to a subsequence if necessary that $\{|\nabla u^{\varepsilon_k}|^p\}_{\varepsilon_k}$ is weakly convergent in $L^1(\Omega)$ thus
\begin{eqnarray}
\lim_{\varepsilon_k\rightarrow 0}\int_D \,|\nabla u^{\varepsilon_k}|^p\,d\x=\int_D \,\int_{{\bf R}^n}|\eta|^p\,d\nu_\x(\eta)\,d\x=\int_D \,\int_Q\,|P(\z)\nabla u^H(\x)|^p\,d\z\,d\x,
\label{optimalityper}
\end{eqnarray}
and optimality follows. 
Last, it follows immediately from Proposition 6.5 of \cite{pedregal}
that for every Caratheodory function $\psi(\x,\eta)$ satisfying the growth condition (\ref{equality})
that (on passage to a subsequence if necessary) 
\begin{eqnarray}
\lim_{\varepsilon_k\rightarrow 0}\int_D\,\psi(\x,\nabla u^{\varepsilon_k})\,d\x=\int_D\,\int_{{\bf R}^n}\psi(\x,\eta)\,d\nu_\x(\eta)\,d\x
\label{cara}
\end{eqnarray}
and (\ref{caratheodory}) follows since (\ref{identify}) implies that
\begin{eqnarray}
\int_D\,\int_{{\bf R}^n}\psi(\x,\eta)\,d\nu_\x(\eta)\,d\x=\int_D\,\int_Q\,\psi(\x,P(\z)\nabla u^H(\x))\,d\z\,d\x.
\label{caratheoident}
\end{eqnarray}

\section{\bf Acknowledgments}
This research effort is sponsored by NSF through grant DMS-0406374 and by the Air
Force Office of Scientific Research, Air Force Material Command USAF, under grant numbers F49620-02-1-0041
and FA9550-05-1-0008.
The US Government is authorized to reproduce and distribute reprints for governmental purposes
notwithstanding any copyright notation thereon. The views and conclusions herein are those of
the authors and should not be interpreted as necessarily representing the official policies or
endorsements, either expressed or implied of the Air Force Office of Scientific Research or the US Government.

\end{document}